# Valuation and Divisibility

Nicholas Phat Nguyen[1]

**Abstract.** In this paper, we explain how some basic facts about valuation can help clarify many questions about divisibility in integral domains.

## I. INTRODUCTION

As most people know from trying to understand and digest a bewildering array of mathematical concepts, the right perspective can make all the difference. A good perspective helps provide a clarifying thread through a maze of seemingly disconnected matters, as well as an economy and unity of thought that makes mathematics much more understandable and enjoyable.

In this note, we want to show how some basic facts about valuation can help provide such a unifying and clarifying perspective for the study of divisibility questions in integral domains. The language and concepts of valuation theory are rather more expressive and intuitive than the purely algebraic questions of divisibility, and through the link with ordered groups, provide some ready-made and effective tools to study and manage questions of divisibility. Most of this is well-known to the specialists, but deserves to be better known by a wider mathematical public.

In what follows, we will focus on an integral domain $A$ whose field of fractions is $K$, unless indicated otherwise. A non-zero element $x$ of $K$ is said to divide another non-zero element $y$ if $y = ax$ for some element $a$ in the ring $A$. It follows from this definition that $x$ divides $y$ if and only if the principal fractional ideal $Ax$ contains the principal fractional ideal $Ay$.

[1] E-mail address: nicholas.pn@gmail.com

Moreover, two elements divide each other if and only if they generate the same fractional ideals, which means they are equal up to a unit factor.

Questions about divisibility in $A$ are then naturally related to the set $D$ of principal fractional ideals $\neq (0)$. Two elements of $K$ that generate the same non-zero fractional ideal are equivalent as far as divisibility is concerned, and are said to be associates of each other. The set $D$ has a natural group structure induced by multiplication of ideals, and whose neutral element is the ideal $A$ itself. The group $D$ is naturally isomorphic to the quotient group $K^*/U$, where $K^*$ is the multiplicative group of non-zero elements in $K$ and $U$ is the group of units in the ring $A$.

The quotient group $K^*/U$ has a partial ordering given by the relation of divisibility. Specifically, the relation ($x$ divides $y$) on $K^*$ is reflective and transitive, and when we pass to the quotient group $K^*/U$ it also becomes anti-symmetric, because two elements of $K^*$ that divide each other represent the same class in $K^*/U$.

In the group $D$, the equivalent order relation is containment or reverse inclusion. For any two non-zero elements $x$ and $y$ in $K$, we have the following equivalent relations in $D$ and $K^*/U$ to express the fact that $x$ divides $y$:

$Ax \supset Ay$ in $D$ is equivalent to ($x$ mod $U$) divides ($y$ mod $U$) in $K^*/U$

We will express these order relations by writing $Ax \leq Ay$ or ($x$ mod $U$) $\leq$ ($y$ mod $U$). For convenience, we sometimes will also write $x \leq y$ to express the relation of divisibility in $K^*$, although it should be noted that we do not really have an ordering until we pass to the quotient group $K^*/U$, in order to avoid the situation of $x \leq y$ and $y \leq x$, but $x \neq y$.

This ordering is compatible with the group structure on $D$ and $K^*/U$, in the sense that the ordering is invariant under composition with any element of the group, i.e., composing both sides of an inequality by any element of the group will not change the inequality. Specifically, if $Ax \leq Ay$, then for any element $z$ in $K^*$, we have $Axz \leq Ayz$.

A *commutative* group with a compatible partial ordering is known as an ordered group. In this case, $D$ and $K^*/U$ are isomorphic ordered groups.

When talking about ordered groups, we will generally use additive notation in order to make things easier to follow, especially as our experience and intuition with ordered groups comes first from the additive groups of $\mathbb{Z}$, $\mathbb{Q}$ and $\mathbb{R}$, with the usual ordering. In particular, we will refer to elements $x \geq 0$ as positive (and $x > 0$ as strictly positive). Likewise for negative and strictly negative numbers. The zero element in the group is both positive and negative.[2]

With the ordering induced by divisibility, the integral principal ideals constitute the set of positive elements in $D$. We will refer to $D$, or equivalently $K^*/U$, as the divisibility group of the domain $A$.

## II. VALUATION AND THE LINK WITH ORDERED GROUP THEORY

Let's consider the natural homomorphism $u$ from $K^*$ to $K^*/U$ that maps an element $x$ to its multiplicative coset $x \bmod U$. In order to highlight the fact that we are looking at $K^*/U$ as an ordered group, we will use additive notation for $K^*/U$.

Because $u$ is a homomorphism, we have $u(xy) = u(x) + u(y)$. What about $u(x+y)$? If an element $z$ divides both $x$ and $y$, then obviously $z$ must divide $(x+y)$. That means, if we translate divisibility relations into the equivalent order relations in $K^*/U$,

If $u(z) \leq u(x)$ and $u(z) \leq u(y)$, then $u(z) \leq u(x+y)$.

If $u(x)$ and $u(y)$ have a greatest lower bound or infimum in $K^*/U$, the above condition is equivalent to $\inf(u(x), u(y)) \leq u(x+y)$.

[2] Many people say "positive" for "$> 0$" and "non-negative" for "$\geq 0$". However, it is often much more convenient to include 0 as a positive number to simplify language.

To understand the divisibility relations in $K^*$ relative to $A$, we need to understand the structure of $K^*/U$ as an ordered group (relative to the divisibility order relation). In this regard, any insight provided by the theory of ordered groups is obviously helpful. More generally, we can consider any homomorphim $u$ from $K^*$ to an ordered group with similar conditions on $u(x + y)$ as noted above.

In most situations, it is enough for us to work with an ordered group in which there is *always* the greatest lower bound or infimum of any two elements. In that case, the least upper bound or supremum of any two elements also exists, and such a group is known as a lattice-ordered group, or lattice group for short.

To give the readers a feel for the structural symmetry of an ordered group, we will prove the following.

*Proposition 1*. *Let G be an ordered group, and x and y are two elements in G. If inf(x, y) exists then sup(x, y) also exists, and vice versa. Moreover,* $\inf(x, y) + \sup(x, y) = x + y$.

Proof. Assume that $\inf(x, y)$ exists. Because $x \mapsto -x$ is an order-reversing automorphism of $G$ as an ordered group, it maps $\inf(x, y)$ to $\sup(-x, -y)$, i.e. $\sup(-x, -y) = -\inf(x, y)$.

Let $a$ be any element of $G$. Because translation by $a$ is an order-preserving automorphism of $G$ as an ordered set, $\sup(a - x, a - y)$ exists and is equal to $a + \sup(-x, -y) = a - \inf(x, y)$. When $a = x + y$, we have $\sup(y, x) = (x + y) - \inf(x, y)$, which implies $\sup(x, y) + \inf(x, y) = x + y$. The proof starting with the assumption of $\sup(x, y)$ is similar. ■

We will refer to any homomorphism $u$ from $K^*$ to such a lattice group that satisfies the above condition on $u(x + y)$ as a valuation. Specifically, a valuation $u$: $K^* \to G$ is a mapping from the multiplicative group $K^*$ to a lattice group $G$ that satisfies the following properties:

(1) $u(ab) = u(a) + u(b)$ (addition in the ordered group $G$);
(2) $u(a + b) \geq \inf(u(a), u(b))$ if $a + b \neq 0$.

Condition (1) means that $u$ is a group homomorphism, and condition (2) expresses the basic requirement that if the value $u(x)$ of an element $x$ in $K^*$ is bounded by $u(a)$ and $u(b)$, e.g., if $x$ divides both $a$ and $b$, then it is also bounded by $u(a+b)$, e.g., $x$ must also divide the sum $a+b$.

We do not require the valuation $u$ to be surjective. However, to avoid the trivial case where $u(x) = 0$ for all $x$, we will assume, unless indicated otherwise, that the valuation is not identically zero, i.e., $u(x) \neq 0$ for some $x$.

For convenience, we also define $u(0)$ to be $\infty$, where the infinity symbol $\infty$ is just an element that we add to the group $G$ with the conventions that (i) $\infty >$ any element in $G$ and (ii) $G + \infty = \infty$. With the addition of the infinity symbol $\infty$, properties (1) and (2) hold for any two elements $a$ and $b$ without exception. Moreover, $\mathcal{R}(u) = \{x \in K \mid u(x) \geq 0\}$ is a sub-ring of $K$ that we call the valuation ring of $u$. The group of multiplicative units of $\mathcal{R}(u)$ is exactly the set of elements with valuation zero.

There is a well-developed and extensive theory of lattice groups, and the concept of valuation as we define above can help us analyze and clarify many questions about divisibility in $K^*$ by reference to lattice group theory. In fact, we hope to show below how the concept of valuation can serve as a kind of Ariadne's thread that helps tie together and clarify quite a few results in commutative algebra.

The definition for a valuation as stated above was given in 1932 by Wolfgang Krull, a great and profound mathematician who created many of the central concepts of modern commutative algebra. In his definition, Krull required the value group $G$ to be totally ordered, and for that reason the term Krull valuation is often used to refer the case when the valuation takes its value in a totally ordered group. However, the definition makes sense for the more general case of a lattice-ordered group, where every two elements has a greatest lower bound. As we will see below, there are many ordered groups that arise naturally and that are not totally ordered although they often will satisfy the lattice condition. For ease of reference in this paper, we use the term valuation for the general case where the value group $G$ is a lattice group, and reserving the term Krull valuation for the case when $G$ is totally ordered.[3]

## III. DISCRETE VALUATIONS

If the value group $G$ is isomorphic to $\mathbb{Z}$ (the additive group of rational integers with its natural ordering) then the valuation u is called a discrete valuation. This is the first and still the most well-known type of valuations.

If we look at the field of meromorphic functions defined at a point in $\mathbb{C}$,[4] the order of a function at that point gives us a discrete valuation, whose valuation ring is just the ring of holomorphic functions defined at that point.

Specifically, assume for convenience that the point we are looking at is the origin (the number zero) in $\mathbb{C}$. We know from complex analysis that any non-zero meromorphic function $f$ defined in a neighborhood of 0 would have a Laurent expansion $f(z) = z^n.(a_0 + a_1z + a_2z^2 + ....)$, where the power series $(a_0 + a_1z + a_2z^2 + ....)$ is absolutely and uniformly convergent in a small disc around the origin with the leading coefficient $a_0 \neq 0$ (and therefore representing a holomorphic function $h$ with $h(0) \neq 0$).

We define ord($f$) at 0 to be the integer $n$, so that if $n$ is strictly negative then $f$ has a pole at 0. We have ord($fg$) = ord($f$) + ord($g$) and ord($f + g$) $\geq$ minimum of ord($f$) and ord($g$). This order function is therefore a discrete valuation on the field of meromorphic functions defined at zero.

Similarly, the order of a formal Laurent series over a field $F$ gives us a discrete valuation on $F((X))$, whose valuation ring is the ring of formal power series $F[[X]]$.

[3] What we call a valuation here is also known as a semi-valuation or demi-valuation (demi is a French term meaning half). Some authors also require a valuation to be surjective, but we do not make that assumption in order to expand the scope of applications.

[4] We say a function is defined at a point if it is defined in some unspecified neighborhood of that point.

In the field $\mathbb{Q}$ of rational numbers, if we look at the exponent of a fixed prime $p$ in the primary decomposition of any rational number, that will give us a discrete valuation as well. The valuation ring of such a discrete valuation is not the integers $\mathbb{Z}$, but a bigger ring consisting of all rational numbers which can be expressed as fractions whose denominators have no prime factor $p$. The study of that valuation ring will naturally lead, through a process of completion, to the ring $\mathbb{Z}_p$ of p-adic integers.

If we look at the valuation ring $\mathcal{R}(u)$ of a discrete valuation $u$, any element $x$ can be expressed as product $x = at^n$, where $a$ is a unit (whose valuation is of course zero), $t$ is any element of valuation 1, which is called a uniformizer, and $n$ is the valuation of $x$. The term uniformizer here is a nod to the classical case of complex function theory, where at any point $s$ in its domain of definition, any holomorphic function can be expressed locally (in a sufficiently small neighborhood of the point $s$) as a product of a holomorphic function that is non-zero at $s$, and an integral power of the local variable $(z - s)$ centered at $s$, also called the local uniformizer.

Accordingly, in a discrete valuation ring, all ideals are principal and are powers of the ideal generated by a uniformizer, which is the only non-zero prime ideal. Anytime we have a discrete valuation $u$ on a field $K$, we can define a norm function $N(x)$ for the field $K$ as follows:

$$N(x) = 0 \text{ if } x = 0; \text{ and } N(x) = \exp(-u(x)) \text{ if } x \neq 0.$$

It is easy to see that this norm function is multiplicative $N(xy) = N(x)N(y)$, and satisfies the condition $N(x + y) \leq \max(N(x), N(y))$, which is a stronger version of the triangle inequality. This inequality is known as the ultrametric inequality.

Under this norm, the valuation ring $u(x) \geq 0$ corresponds to the unit disc $N(x) \leq 1$. Moreover, we can define a metric $d(x, y) = N(x - y)$ on the field $K$. The operations of addition, multiplication and inversion $x \mapsto 1/x$ are continuous with respect to this metric topology, and therefore $K$ endowed with this topology is a topological field. So a discrete valuation allows us to apply topology and function theory.

There are many different characterizations of a discrete valuation ring in the literature (e.g., as a local principal ideal domain, a Noetherian domain of dimension 1 which is integrally closed, etc.), and discussions of discrete valuation rings can be found in a wide range of subjects, from commutative algebra, algebraic number theory, algebraic geometry to p-adic analysis and complex function theory. This is a testimony to how important and useful the concept of a discrete valuation has become in many mathematical contexts.

In 1882, Richard Dedekind and Heinrich Weber introduced the concept of discrete valuations in what many people regard as perhaps the greatest mathematical paper of all time.[5] Their paper established what in modern mathematical parlance would be called an equivalence between the category of compact Riemann surfaces and the category of algebraic function fields over the complex numbers, i.e., finite extensions of the field $\mathbb{C}(X)$ of rational functions.

Given a compact Riemann surface $S$, the field of meromorphic functions on $S$ is such an algebraic function field. Conversely, given an algebraic function field over $\mathbb{C}$, how do we construct something equivalent to a Riemann surface $S$? What Dedekind and Weber observed was that each point of such a Riemann surface would give us a discrete valuation on the algebraic function field by taking orders of functions. Their investigations showed that all the essential properties of the Riemann surface are reflected in the discrete valuations of the algebraic number field. Their work is the beginning of modern algebraic geometry and has opened up a vast new landscape in mathematics, where the many geometric ideas from the study of complex functions and complex algebraic curves find new applications in algebra and number theory through the use of valuations.

## IV. EXAMPLES OF MORE GENERAL VALUATIONS

Aside from $\mathbb{Z}$, $\mathbb{Q}$, $\mathbb{R}$ or more generally any totally ordered group, other natural examples of lattice groups include a direct sum of totally ordered groups $G = \oplus_i G_i$, where the ordering on $G$ is defined component-wise, with $(a_i) \geq (b_i)$ if and only if $a_j \geq b_i$ for any index $i$, and similarly for any product $G = \prod_i G_i$ of totally ordered groups.

[5] See *Theory of Algebraic Functions of One Variable*, annotated translation with introduction by John Stillwell of original 1882 paper by Richard Dedekind and Heinrich Weber, American Mathematical Society (2012).

We will refer to a valuation $v$: $K^* \rightarrow G = \bigoplus_i G_i$ as a divisor valuation. If such a divisor valuation has all component groups $G_i$ isomorphic to a totally ordered group $H$, then we say that $v$ is an $H$-divisor valuation. Similarly, we will refer to a valuation $v$: $K^* \rightarrow G = \prod_i G_i$ as a product valuation, or an $H$-product valuation if all component groups $G_i$ are isomorphic to a totally ordered group $H$.

An important foundational result in the theory of lattice groups is that any (commutative) lattice group can be embedded as a sub-lattice and subgroup in a product of totally ordered groups. See, e.g., [Yakabe 1963], section 8. Consequently, a product valuation is the most general example of a valuation.

Below are some examples of valuations.

Unique factorization domains: For any unique factorization domain $A$ (also known as factorial ring), with field of fractions $K$, there is a natural $\mathbb{Z}$-divisor valuation $v$ such that $\mathcal{R}(v) = A$. If $(\pi_i)$ is a complete set of representatives of the irreducible elements in $K$, then the $\mathbb{Z}$-divisor valuation $v$: $K^* \rightarrow \bigoplus_i \mathbb{Z}$ is defined as follows. For any non-zero element $x$ in $K$, we take the i-component of $v(x)$ to be the exponent of the irreducible element $\pi_i$ in the prime factorization of $x$. Note that the valuation $v$ is surjective in this case, i.e., the image $v(K^*)$ of the valuation is the whole group $\bigoplus_i \mathbb{Z}$. In fact, it is easy to see that the ring $A$ is a unique factorization domain if and only if its divisibility group is isomorphic to a $\mathbb{Z}$-divisor group $\bigoplus_i \mathbb{Z}$, or equivalently, that there is *surjective* $\mathbb{Z}$-divisor valuation $v$: $K^* \rightarrow \bigoplus_i \mathbb{Z}$ such that $\mathcal{R}(v) = A$.

Dedekind domains: Consider any Dedekind domain $A$ with field of fractions $K$. For example, $A$ could be the ring of integers in a number field, or the coordinate ring of a regular affine curve. A Dedekind domain has the defining property that any fractional ideal is invertible. There are also many equivalent characterizations. The property that most people associate with Dedekind domains is the following (which can either be regarded as a defining property or as a consequence of some simpler defining conditions): any ideal in a Dedekind domain can be expressed uniquely as a product of prime ideals. This unique factorization is similar to the unique factorization in the ring of integers, except that it is for ideals rather than numbers. This unique factorization for ideals was discovered (and the concept of ideal created in the process) by Richard Dedekind in his landmark study of algebraic numbers.

The primary decomposition of ideals in a Dedekind domain $A$ gives rise to a natural $\mathbb{Z}$-divisor valuation $v$ such that $\mathcal{R}(v) = A$. Specifically, let $(\wp_i)$ be the family of all maximal ideals in $A$. We can define a $\mathbb{Z}$-divisor valuation $v$: $K^* \to \bigoplus_i \mathbb{Z}$ as follows: for any non-zero integral element $x$, the i-component of $v(x)$ is simply the exponent of the ideal $\wp_i$ in the primary ideal decomposition of the principal ideal $Ax$. We can extend that definition to all non-zero elements of the field in a natural way.

Prüfer domains: If an integral domain $A$ has the property that any finitely-generated fractional ideal is invertible, we have the following natural valuation defined on the field $K$ of fractions of $A$. Such an integral domain is known as a Prüfer domain and can be thought of as a generalization of Dedekind domain, with Dedekind = Prüfer + Noetherian condition. Let $G$ be the set of all finitely-generated fractional ideals of $A$. $G$ is a commutative group under multiplication of ideals, with $A$ as the neutral or identity element. We can order the elements of $G$ by reverse inclusion, i.e., we say $I \leq J$ if $I \supset J$. Such an ordering is compatible with ideal multiplication. Moreover, given any finitely generated fractional ideals $I$ and $J$, the ideal $I + J$ is the infimum of $I$ and $J$ in the ordered group $G$. Hence $G$ is a lattice group. Define $v$: $K^* \to G$ by letting $v$(x) be the ideal $Ax$. $v$ is a valuation.

The above example of valuation on a Prüfer domain can be regarded as a product valuation. For each maximal ideal $\wp_i$ in the domain $A$, the local ring at $\wp_I$ has the same property that any finitely generated fractional ideal is invertible, which means such a finitely generated ideal must be principal because any invertible fractional ideal of a local integral domain is principal. Accordingly, each local ring of $A$ at a maximal ideal $\wp_i$ is a Krull valuation ring, with the value group being the group $G_i$ of principal fractional ideals ordered by reverse inclusion. We have then an order-preserving group homomorphism from the ordered group $G$ of finitely generated fractional ideals in $A$ to the product $\prod_i G_i$, where each ideal $I$ is mapped to its vector of localizations at all the maximal ideals of $A$. If we compose the natural valuation from $K$ to $G$ defined above with this homomorphism, we have a valuation from $K$ into $\prod_i G_i$.

If the Prüfer domain $A$ is also Noetherian and hence a Dedekind domain, each of the component value group $G_i$ is isomorphic to $\mathbb{Z}$. Moreover, each principal ideal $Ax$ is contained in only a finite number of maximal ideals of $A$. Hence the natural valuation defined by finitely-generated fractional ideals is actually a $\mathbb{Z}$-divisor valuation in the case of a Dedekind domain.

<u>Ring of algebraic integers</u>: Let $\overline{\mathbf{Q}}$ be an algebraic closure of $\mathbb{Q}$, e.g., the set of all algebraic numbers in the field $\mathbb{C}$ of complex numbers. The ring of all algebraic integers has the Bézout property that any finitely generated ideal in that ring is principal. Because principal ideals are invertible, the ring of all algebraic integers is therefore a Prüfer domain. The natural product valuation on $\overline{\mathbf{Q}}$ (as the field of fractions of a Prüfer domain) can be defined explicitly as follows.

For each maximal ideal $\wp_i$ in the ring of all algebraic integers, there is a Krull valuation $v_i$ taking value in the ordered group $\mathbb{Q}$. Specifically, for each non-zero algebraic number $x$, consider the maximal ideal $\wp_i \cap K$ in any number field $K$ containing $x$, such as $K = \mathbb{Q}(x)$. That maximal ideal defines a discrete valuation $u$ on $K$. We define $v_i(x) = u(x)/e$, where $e$ is the ramification index of $\wp_i \cap K$ over $\mathbb{Z}$. It follows from the transitivity of ramification indices for number fields that the above formula gives the same valuation for any number field $K$ containing $x$, and therefore is well-defined. Let $v(x) = \prod_i v_i(x)$, where the product is taken over the index set of all maximal ideals $\wp_i$ in the ring of all algebraic integers. $v$ defines a **v**aluation on $\overline{\mathbf{Q}}$ such that $\mathcal{R}(v)$ is the ring of all algebraic integers.

<u>Meromorphic functions on a Riemann surface</u>: Let $\mathcal{M}(X)$ be the field of meromorphic functions on a Riemann surface $X$. For any point $p$ of $X$, there is a natural discrete valuation $v_p$ defined by taking the order of a meromorphic function at $p$. There is therefore a natural $\mathbb{Z}$-product valuation on $\mathcal{M}(X)$ where each component valuation is the order valuation at each point of the Riemann surface. The valuation ring of this $\mathbb{Z}$ -product valuation is the subring A of holomorphic functions.

If the Riemann surface is compact, then the above $\mathbb{Z}$-product valuation is actually a $\mathbb{Z}$-divisor valuation $v$: $K^* \to \oplus_p \mathbb{Z}$ because each meromorphic function has a finite number of poles and zeros on $X$. This divisor valuation is not surjective because for any function $f$ we have $\sum v_p(f) = 0$ (any meromorphic function on a compact Riemann surface has as many zeros as poles, counting multiplicities). The valuation ring in this case, i.e. the ring of holomorphic functions on $X$, consists of just the constant functions.

## V. UNIQUE FACTORIZATION DOMAINS

We outline here a dictionary between divisibility and valuation concepts that would be helpful in recognizing the key features of unique factorization in their many guises.

An element $z$ of $A$ is regarded as irreducible if it has no proper divisors. That is equivalent to the ideal $Az$ being maximal (with respect to set inclusion) relative to all proper integral principal ideals, or that $Az$ is minimal in the set of strictly positive elements of $D$. In modern parlance, $z$ is also called an atom if it is irreducible.

Factorization of an element of $A$ into irreducibles is equivalent to decomposition of a positive element as a finite sum of minimal elements in the ordered set $D$ or $K^*/U$, which can easily be seen as equivalent to the Artinian condition on the ordered set, namely any decreasing sequence of positive elements must become stationary at some point, or that any non-empty set of positive elements must have a minimal member (no other element in the set being strictly smaller). [6]

---

[6] In terms of inclusion in the set integral principal ideals, that means we must have the Noetherian condition, namely any increasing sequence of integral principal ideals becomes stationary at some point, or equivalently that any non-empty family of such ideals must have a maximal element.

We will refer to any integral domain that satisfies the above equivalent conditions for factorization as a factorization domain.

Factorization may be impossible if the ring *A* is really big. For example, if *A* is the valuation ring of a surjective valuation $u$: $K^* \to \mathbb{Q}$, then it is easy to see that there can be no irreducible element because the ordered group $\mathbb{Q}$ has no strictly positive minimal element, so factorization does not even get off to a start. Similarly, the total ring of all algebraic integers has no irreducible element because given any algebraic integer *x* always has proper divisors, such as ($x^{1/2}$). However, for Noetherian integral domains, which is the type of integral domains we encounter in most applications, factorization is automatic. Also, for rings such as $K[X_i]$ (polynomial ring in arbitrary number of variables), which is a direct limit of smaller rings stable under factorization, the condition is also satisfied if it is satisfied in each ring of the directed family.

Uniqueness of factorization depends, as has been known since the time of Euclid, on irreducible elements also having the prime divisor property. Specifically, an element *p* has the prime divisor property if whenever *p* divides a product *ab*, *p* must either divide *a* or *b*. This prime divisor property is equivalent to the ideal *Ap* being a prime ideal, i.e. that *p* is a prime element. It is easy to see that a prime element *p* is irreducible. The converse is equivalent to uniqueness of factorization and is often a more delicate condition to verify than factorization.

Here the theory of ordered groups can be of significant help. In the language of ordered groups, the prime divisor property is equivalent to the following: if a minimal element *p* in the set of strictly positive elements is bounded by a sum $(a + b)$ of two positive elements, then *p* must be bounded by either *a* or *b*. It happens that if the ordered group has the lattice property, then this prime divisor property is automatically satisfied for such a minimal element *p*, if it exists.[7]

---

[7] A lattice group may not have a minimal element, e.g., the rational numbers $\mathbb{Q}$ in the usual ordering, but if it has a minimal element, then the minimal element must necessarily have the prime divisor property.

Indeed, a basic decomposition theorem for lattice groups says that if $p \leq a + b$ ($p$, $a$, $b$ being positive elements) then we must have $p = c + d$ ($c$ and $d$ positive) with $c \leq a$ and $d \leq b$. See [Bourbaki A2], chapter 6, no. 10. If $p$ is irreducible, that means either $c$ or $d$ must be zero, which means we have $p = c$ or $p = d$, and therefore $p \leq a$ or $p \leq b$.

Because the divisibility group $D$ of a unique factorization domain is isomorphic (as an ordered group) to a $\mathbb{Z}$-divisor group $\oplus_i \mathbb{Z}$, which is of course a lattice group, uniqueness of factorization is then equivalent to $D$ having the lattice condition. That means for any two elements $x$, $y$ in K*, we can always find their greatest lower bound inf($x$, $y$) or least upper bound sup($x$, $y$).

Recall that the greatest lower bound $a$ of $x$ and $y$ is defined by the property that $a$ is a common lower bound of both $x$ and $y$, and for any common lower bound $c$ of $x$ and $y$, we have $c \leq a$. So for two elements $x$, $y$ in $K^*$, their greatest lower bound inf($x$, $y$) is the element $a$ of K* (determined up to a unit factor) such that $a$ divides both $x$ and $y$, and any common divisor $c$ of $x$ and $y$ also divides $a$.

Here the reader may be tempted to say, oh, but that is equivalent to $a$ being the greatest common divisor of $x$ and $y$ when $x$ and $y$ are integral elements. But there is a subtlety here that we must think through. For integral elements $x$ and $y$, the greatest common divisor of $x$ and $y$ is defined (up to a unit factor) by reference to common divisors that are integral elements. In other words, gcd($x$, $y$) = infimum of $x$ and $y$ in the set of integral elements, i.e., positive elements of $K$ under the divisibility order relation. If we consider all elements of $K$, $\inf_K(x, y)$ may not exist.

However, it turns out that if $\gcd(x, y) = \inf_A(x, y)$ exists for *any* two elements $x$, $y$ of $A$, then $\inf_K(x, y)$ also exists and is equal to $\inf_A(x, y)$. The reason is that any element of $K^*$ is a fraction of two integral elements, so if $c$ is a lower bound of both $x$ and $y$, we can change $c$ into an element of $A$ with suitable multiplication by an integral element $h$. Using additive notation, we have $c + h \leq \inf_A(x + h, y + h)$, which is equal to $\inf_A(x, y) + h$, where $\inf_A$ is taken in the set of integral elements and assumed to exist by hypothesis. Reverse translation by $h$ gives us $c \leq \inf_A(x, y)$, so $\inf_A(x, y)$ is also the infimum $\inf_K(x, y)$ in K*.

For $\sup(x, y)$, the situation is simpler, since any common upper bound of $x$ and $y$ must necessarily be positive, and so $\sup_A(x, y)$ taken in the set of positive elements is the same as $\sup_K(x, y)$ taken in the entire ordered set.

A little reflection shows that the least upper bound $\sup(x, y) = \text{lcm}(x, y)$, if it exists, must correspond to the intersection $Ax \cap Ay$. Specifically, let $\sup(x, y) = d$, then $x \leq d$ and $y \leq 0$ means $Ad \subset Ax$ and $Ad \subset Ay$ or $Ad \subset Ax \cap Ay$. At the same time, any element $t$ in $Ax \cap Ay$ is an upper bound of both $x$ and $y$, so we also have $d \leq t$, meaning $At \subset Ad$, so that $Ax \cap Ay$ is also contained in $Ad$.

On the other hand, all we can say about the greatest lower bound $\inf(x, y) = \gcd(x, y)$, if it exists, is that it must correspond to the intersection of all principal ideals containing the ideal generated by $x$ and $y$. If the integral domain satisfies the Bézout property that any finitely-generated ideal is principal, then $\inf(x, y)$ is indeed the same as the principal ideal generated by $x$ and $y$.

If $A$ is the valuation ring $\mathcal{R}(v)$ for a surjective valuation $v$: $K^* \to G$ then for any elements $x$, $y$ in $A$, $\gcd(x, y)$ exists because it must be the element (up to a unit factor) corresponding the valuation $\inf(v(x), v(y))$. Similarly, $\text{lcm}(x, y)$ exists because it must be the element (up to a unit factor) corresponding the valuation $\sup(v(x), v(y))$.

To summarize, uniqueness of factorization is equivalent to any one of the following conditions, most of which derive from properties of the divisibility group viewed as an ordered group:

- Any principal ideal generated by an irreducible element is a prime ideal. (If A is a Noetherian domain, it follows from Krull's principal ideal theorem that this is equivalent to any prime ideal of height 1 being principal.)
- The divisibility group $D$ or $K^*/U$ is a lattice group (relative to the order of divisibility).
- $A$ is the valuation ring for a surjective valuation $v$: $K^* \to G$.
- For any two elements $x$ and $y$ of $A$, $\inf_A(x, y) = \gcd(x, y)$ exists.

- For any two elements $x$ and $y$ of $A$, $\sup_A(x, y) = \operatorname{lcm}(x, y)$ exists.
- For any two elements $x$ and $y$ of $A$, the intersection $Ax \cap Ay$ is a principal ideal.
- For any two elements $x$ and $y$ of $A$, the intersection of all principal ideals containing both $x$ and $y$ is a principal ideal.

A particularly simple but very useful situation is when a factorization domain has essentially only one irreducible element.

*Proposition 2*. *Let A be a factorization domain. The following conditions are equivalent.*

*(a) A has only one irreducible element (up to a unit factor).*
*(b) A is a discrete valuation ring.*
(c) *If A has only one maximal ideal* $\mathbf{m}$*, and* $\mathbf{m}$ *is a principal ideal.*

Proof. *(a) implies (b)*: Expressed in terms of the divisibility group of $A$, this means if an ordered group has a unique (strictly positive) minimal element and if every element in the ordered group is generated by minimal elements, then the ordered group is isomorphic to $\mathbb{Z}$. Indeed, the cyclic subgroup generated by this minimal element must be the whole group, and it is torsion-free because we are in an ordered group.

*(b) implies (c)*: If $t$ is a uniformizer of the discrete valuation ring $A$, then the ideal $(t)$ is the only maximal ideal of $A$.

*(c) implies (a)*: Let $\mathbf{m} = (t)$. Because $\mathbf{m}$ is maximal, the element $t$ is irreducible and prime. Moreover, any irreducible element is divisible by $t$ because the ideal $(t) = \mathbf{m}$ is the only maximal ideal, and so any irreducible element must be associated to $t$. ■

*Corollary 2.1* *Let A be a factorization domain. The localization of A at a prime element is a discrete valuation ring.*

Proof. Recall that the localization of an integral domain $A$ at a prime element $p$ is the ring of all fractions with denominators not divisible by $p$. (Such a ring is well-defined because $p$ has the prime divisor property.) That ring of fractions is a local ring whose maximal ideal is the principal ideal generated by $p$. Any element in this local ring at $p$ can be expressed as power of $p$ multiplied by a unit, and hence it is a unique factorization domain. The result now follows from Proposition 2. ■

We note below some necessary properties of a factorization domain with a finite number of irreducible elements. In such a domain, it is clear that any non-zero prime ideal must be generated by a finite number of irreducible elements. Therefore the domain must also be Noetherian, by a well-known theorem of P.M. Cohn.

*Proposition 3*. *Let A be a Noetherian domain and K the field of fractions of A. If A has a finite number of irreducible elements, then A has the following properties.*

*(a) The field K is a finite ring extension of A (i.e., K is generated as a ring over A by a finite number of elements).*
*(b) The intersection of all non-zero prime ideals of A is a non-zero ideal.*
*(c) A has a finite number of non-zero prime ideals, and all of them are maximal.*

Proof. It is easy to see that if $A$ has a finite number of irreducible elements, then the field $K$ is generated by adjoining the element $1/d$ to $A$, where $d$ is the product of all the irreducible elements in $A$.

The conditions of (a), (b) and (c) are actually equivalent and follow from the characterizations by Emil Artin and John Tate of a Noetherian domain A whose field of fractions is a finite ring extension of A. See Theorem 4 of [Artin & Tate 1951]. ■

In particular, a factorization domain whose Jacobson radical (intersection of all maximal ideals) is zero must have an infinite number of irreducible elements. This is essentially equivalent to what Professor Pete L. Clark called the Euclidean criterion for irreducibles (as applied to a factorization domain) in a 2017 article in the American Mathematical Monthly. See [Clark 2017].

A factorization domain with a finite number of irreducible elements is known as a Cohen-Kaplansky domain. These domains were completely characterized in 1992 after being largely forgotten for many years. For more detail, please see [Clark 2017] and the references cited in that article.

## VI. THINKING IN TERMS OF VALUATION

Modern algebra is dominated by the conceptual approach pioneered by Richard Dedekind. Under this approach, the key concept in the study of commutative rings is the concept of ideals. There is also another approach pioneered by Leopold Kronecker, a contemporary of Dedekind, who advocated the use of algorithms and constructive methods. The constructive methods of Kronecker have been out of sight for a long time, although they are beginning to be better understood and appreciated in recent times thanks in part to the tireless effort by Professor Harold Edwards to explain the constructive ideas and methods of Kronecker through a series of beautiful books and articles.

The use of valuations can be thought of as an intermediate approach between the conceptual ideal-centric approach of Dedekind and the constructive formula-centric approach of Kronecker. Each valuation on a field is a kind of yardstick that allows us to measure and compare different elements, which helps clarify and simplify many issues of divisibility.

Perhaps the most concise way to describe a unique factorization domain is as the valuation ring of a surjective $\mathbb{Z}$-divisor valuation. If the $\mathbb{Z}$-divisor valuation is not necessarily surjective, then we have what is known as a Krull domain. This view of Krull domain as the valuation ring of a $\mathbb{Z}$-divisor valuation is also perhaps the most concise way to describe this important class of integral domains.

We can define the class group of a valuation $u$: $K^* \to G$ as the quotient group $G/\mathrm{u}(K^*)$. A unique factorization domain is then a $\mathbb{Z}$-divisor valuation ring with trivial class group. The class group of the natural $\mathbb{Z}$-divisor valuation on a Dedekind domain is the same as the ideal class group defined taking the group of all fractional ideals modulo the group of principal ideals.

The integral closure of a subring $B$ of a field $K$ can also be characterized by valuation rings as follows. An element $x$ of $K$ is integral over $B$ if and only if for any Krull valuation $u$ on $K$ with $u(B) \geq 0$, we also have $u(x) \geq 0$. This characterization leads immediately to the following, which we note here for future reference.

*<u>Proposition 4</u>. Let B be a subring of a field K.*

*(a) B is a valuation ring of K if and only if it is integrally closed in K.*
*(b) For any elements x and y in K*, x divides y relative to the integral closure of B in K if and only if $w(x) \leq w(y)$ for all valuations w with $w(B) \geq 0$.*

<u>Proof</u>.

(a) Let $u(x) = (w_i(x))$ where the valuations $w_i$ runs through all the Krull valuations of $K$ that take positive values on $B$. It is clear that $u$ is product valuation whose valuation ring $\mathcal{R}(u)$ is the intersection of all the Krull valuation rings $\mathcal{R}(w_i)$. By the above characterization of integral elements over $B$, $\mathcal{R}(u)$ is the integral closure of $B$ in $K$. So if $B$ is integrally closed in K, it is the valuation ring $\mathcal{R}(u)$.

It is easy to see that each Krull valuation ring $\mathcal{R}(w_i)$ is integrally closed in $K$ by applying the Krull valuation $w_i$ to any integral dependence equation. It follows that a valuation ring $\mathcal{R}(u)$ is also integrally closed in $K$ because it is the intersection of all the rings $\mathcal{R}(w_i)$.

(b) Let $z = y/x$. The condition $w(x) \leq w(y)$ is equivalent to $w(z) \geq 0$. If $w(z) \geq 0$ for all valuations $w$ with $w(B) \geq 0$, in particular for all Krull valuations $w$, then $z$ must be integral over $B$.

Conversely, if $z$ is integral over $B$, then the equation of integral dependence forces any Krull valuation $w(\text{z})$ to be $\geq 0$ if $w(B) \geq 0$. Any valuation $u$ can be expressed as a product valuation whose components are Krull valuations, so $u(z)$ must also be $\geq 0$ if $u(B) \geq 0$. ■

Thinking in terms of valuation can make quite a few things more transparent. For example, consider the following very useful theorem by Masayoshi Nagata.

*Theorem 5 (Nagata)*: *Let A be a factorization domain, and let S be a multiplicative subset of A generated by prime elements. If the ring of fractions $S^{-1}A$ is a unique factorization domain, then A is also a unique factorization domain.*

Proof: For each prime element $p$ in $S$, the ring $A$ localized at $p$ is a discrete valuation ring by Corollary 2.1. So each such prime element $p$ gives us a discrete valuation.

Because $S^{-1}A$ is unique factorization domain, it is the valuation ring of a surjective $\mathbb{Z}$ -divisor valuation. If we put this $\mathbb{Z}$-divisor valuation together with all the discrete valuations defined by prime elements of $S$, we have a larger $\mathbb{Z}$-divisor valuation whose valuation ring is just $A$ itself. Moreover, this larger $\mathbb{Z}$-divisor valuation is also surjective because each prime generator of $S$ has value 1 at the discrete valuation associated to $p$, but has value zero at the $\mathbb{Z}$-divisor valuation of $S^{-1}A$ and at all the discrete valuations associated to the other prime generators of $S$. So $A$ is also a unique factorization domain. ■

We now consider principal ideals in a valuation ring. This will be handy later.

A discrete valuation ring is a principal ideal domain. Any finitely-generated ideal in a Krull valuation ring is a principal ideal. For a more general valuation ring, we need an additional condition to assure that a finitely-generated ideal is principal. Specifically, we need the following Bézout property for the valuation $v$ (in analogy with the familiar Bézout identity for the ring of rational integers or any principal ideal domain):

*(Bézout property) for any elements x and y in the valuation ring $\mathcal{R}(v)$, $\inf(v(x), v(y)) = v(cx + dy)$ for some coefficients c and d in $\mathcal{R}(v)$.*

*Proposition 6*. *Let A be an integral domain contained in a field K (we do not assume that K is the field of fractions of A). Assume $A = \mathcal{R}(u)$ for some valuation u on K.*

(a) *Any finitely-generated ideal of A is a principal ideal if and only if the valuation u satisfies the Bézout property. In other words, A is a Bézout domain if and only if its valuation satisfies the Bézout property.*

(b) *If u is a $\mathbb{Z}$-divisor valuation, then A is a principal ideal domain if and only if the valuation u satisfies the Bézout property.*

Proof of (a). Assume that the valuation u satisfies the Bézout property. It is sufficient to show that any ideal $J$ of $A = \mathcal{R}(u)$ that is generated by two elements must be principal.

Let $x$ and $y$ be the generating elements of $J$. By the Bézout property, there is an element $m$ in $J$ such that $u(m) = \inf(u(x), u(y))$. Since $J$ is generated by $x$ and $y$, any element $z$ of $J$ has the form $cx + dy$ with some coefficients $c$ and $d$ in $A$, and therefore the valuation $u(z) \geq \inf(u(cx), u(dy)) \geq \inf(u(x), u(y)) = u(m)$. Because $z/m$ has positive valuation and hence belongs to $A$, $z = (z/m).m$ is in the principal ideal generated by $m$. Therefore $J$ is a principal ideal generated by $m$.

Conversely, assume that any finitely-generated ideal of $A$ is principal. We need to show that $A$ satisfies the Bézout property. Let $J$ be the ideal generated by $x$ and $y$. By assumption, $J$ is generated by a single element $m$. Because $x$ and $y$ each is a product of $m$ by an element of positive valuation, we certainly have $u(m) \leq u(x)$ and $u(m) \leq u(y)$, and therefore $u(m) \leq \inf(u(x), u(y))$. On the other hand, $m$ is in the form $cx + dy$ for some coefficients $c$ and $d$ in $A$, and so $u(m) \geq \inf(u(cx), u(dy)) \geq \inf(u(x), u(y))$. Accordingly, $u(m) = \inf(u(x), u(y))$.

Proof of (b). If $u$ is a $\mathbb{Z}$-divisor valuation $u$: $K^* \rightarrow \oplus_i \mathbb{Z}$, observe that any non-empty subset of the positive cone of $\oplus_i \mathbb{Z}$ (i.e., the set of elements $\geq 0$) must have a minimal element. Indeed, for any element α in such a subset, there are only a finite number of elements in the positive cone of $\oplus_i \mathbb{Z}$ that are strictly smaller than α. Hence by a finite process of descent, we can always arrive at a minimal element in the given subset (i.e., with nothing that is strictly smaller).

If $A$ is a principal ideal domain then by proposition (a) above, $A$ certainly satisfies the Bézout property. Conversely, assume that $A$ satisfies the Bézout property. Let $J$ be any ideal of $A$. As noted, we can choose an element $m$ of $J$ such that the value $u(m)$ is minimal in the set $u(J)$. Let $z$ be any element of $J$. By the Bézout property, there is an element $t$ in $J$ such that $u(t) = \inf(u(m), u(z))$. Because $u(m)$ is a minimal value in the set $u(J)$, we must have $u(t) = u(m) \leq u(z)$. Accordingly, $z$ must be in the principal ideal generated by $m$. Since $z$ is arbitrary, $J$ must be generated by $m$. ■

## VII. EXTENSION OF VALUATION TO THE FIELD OF RATIONAL FUNCTIONS

The Kronecker constructive approach generally requires working with multivariable polynomials. A question in a base field $K$ can often be solved by first constructing a polynomial in one or several variables over $K$, and then using calculations with that polynomial in a function field over $K$ to lead us to an answer in the base field. In many cases, such constructive approach corresponds to extending a valuation from a base field $K$ to a function field over $K$.

We consider an integral domain $A$ contained in the field $K$ (we do not assume that $K$ is the field of fractions of $A$). Assume that we have a valuation $v$ on $K$ such that $A = \mathcal{R}(v)$. We want to know how to extend the valuation $v$ on $K$ to a valuation on the function field $K(X)$.

If we look at $K(X)$ as the field of fractions of the principal ideal domain $K[X]$, there is a natural $\mathbb{Z}$-divisor valuation $t$ on $K(X)$ defined by factorization into irreducible polynomials. For that valuation $t$, $t(z) = 0$ for any non-zero element $z$ of $K$, so the valuation $t$ is not an extension of $v$.

<u>*Theorem 7*</u>. *Any valuation v on K can be extended to a valuation w on the field of rational functions K(X) such that for any polynomial* $f(X) = a_mX^m + .... + a_1X + a_0$*, we have*

$$w(f) = \inf(v(a_m), ..., v(a_1), v(a_0)).$$ [8]

<u>Proof</u>. The function $w$ as defined is obviously an extension of $v$ in the sense that we have $w(c) = v(c)$ for any element $c$ in $K$ (regarded as a constant polynomial in $K[X]$).

We will need to show that $w$ as defined above for the polynomial ring $K[X]$ has the multiplicative property, namely $w(fg) = w(f) + w(g)$ for any two polynomials $f$ and $g$. If that holds, then it is easily seen that $w$ can be extended to the field of rational functions $K(X)$ with the same multiplicative property.

[8] Rational functions in the field $K(X)$ are also known as rational fractions with coefficients in $K$.

By construction, we have $w(f + g) \geq \inf(w(f), w(g))$ for any two polynomials $f$ and $g$. If the multiplicative property holds, then the same inequality can easily be extended to any two rational functions $f$ and $g$ in the field $K(X)$ if we express $f$ and $g$ with the same denominator.

Therefore, the theorem is established if we can prove the following lemma:

*<u>Lemma (Gauss - Kronecker)</u>: $w(fg) = w(f) + w(g)$ for any two polynomials f and g.*

This lemma is in fact a generalization of the well-known lemma due to Gauss that the content of a product of two polynomials in $\mathbb{Z}[X]$ is equal to the product of the contents of each individual polynomial. Gauss's lemma readily generalizes to the case of any unique factorization domain. In his ideal theory, Kronecker also generalized Gauss's lemma to the case of any Dedekind domain, where the content of a polynomial is no longer just an element of the ground field (up to associates), but a fractional ideal. See [Flanders 1960]. For any Dedekind domain $A$ with field of fractions $K$, let $v$ be the natural $\mathbb{Z}$-divisor valuation on $K$ taking values in the group of fractional ideals of $A$. Then the content of a polynomial $f(X)$ in the sense of Gauss or Kronecker is essentially the same as the value $w(f)$ defined above.

It is enough to show that the multiplicative property holds when $v$ is a Krull valuation, since the multiplicative property would then extend component-wise to the more general case of a product valuation, which means any valuation (because any lattice group can be embedded as a subgroup and sub-lattice in a product of totally-ordered groups).

In the case of a Krull valuation, because the value group $G$ is totally ordered, one of the coefficients of the polynomial $f = a_m X^m + \dots + a_1 X + a_0$, say $a_i$, must have the minimum value $\inf(v(a_m), \dots, v(a_1), v(a_0)) = w(f)$. Similarly, let $b_j$ be the coefficient of $g = b_n X^n + \dots + b_1 X + b_0$ with minimum value $= w(g)$. We can also choose the indices $i$ and $j$ so that each is the largest index with that minimal property. Then in the polynomial $fg = c_p X^p + \dots + c_1 X + c_0$, the coefficient $c_{i+j}$ clearly has the value $w(f) + w(g)$. Moreover, all the other coefficients of $fg$ have values that are at least $w(f) + w(g)$. Accordingly, $w(fg) = v(c_{i+j}) = w(f) + w(g)$. ■

The polynomial $X$ obviously has the valuation $w(X) = 0$. Accordingly, $X$ and its powers are units in the valuation ring $\mathcal{R}(w)$.

NOTE: Although the theorem is stated for the case of a function field in one variable, it can be readily extended to the case of an arbitrary family of variables.

<u>Corollary 7.1</u>. *The valuation w defined on K(X) satisfies the Bézout property.*

<u>Proof</u>. For any two elements $f$ and $g$ in $K(X)$, we can always express them with the same denominator. Let $f = p/r$ and $g = q/r$, where $p$, $q$ and $r$ are polynomials in $K[X]$. We have $\inf(w(f), w(g)) = \inf(w(p) - w(r), w(q) - w(r)) = \inf(w(p), w(q)) - w(r)$. Hence, it is enough to show that there are polynomials $c$ and $d$ in $\mathcal{R}$(w) such that

$$w(cp + dq) = \inf(w(p), w(q)).$$

Take $c = X^h$, where $h$ is any integer > deg($q$), and $d = 1$. Then all the coefficients of $p$ and $q$ appear as separate coefficients in the polynomial $cp + dq = X^h p + q$. In that case, we clearly have $w(cp + dq) = \inf(w(p), w(q))$. ■

<u>Corollary 7.2</u>. *Any finitely generated ideal of $\mathcal{R}(w)$ is principal, i.e., $\mathcal{R}(w)$ is a Bézout domain.*

This is immediate from Corollary 6.1 and Proposition 5(a). ■

For example, let $v$ be the valuation on $\overline{\mathbb{Q}}$ that we describe earlier such that $\mathcal{R}(v)$ is the ring of all algebraic integers. If the valuation $w$ on $\overline{\mathbb{Q}}$(X) extends $v$, then in the valuation ring $\mathcal{R}(w)$, every finitely generated ideal is principal.

Similarly, putting Corollary 6.1 and Proposition 5(b) together, we have:

*Corollary 7.3 (Kronecker).* *If $v$ (and hence $w$) is a $\mathbb{Z}$-divisor valuation, then $\mathcal{R}(w)$ is a principal ideal domain.*

For example, let $A$ be a Dedekind domain with field of fractions $K$, and let $v$ be the natural $\mathbb{Z}$-divisor valuation defined by primary ideal decomposition. Then the subring $\mathcal{R}(w)$ of the field of rational functions $K(X)$ is a principal ideal domain. This is a key result of the classical ring theory developed by Kronecker. See [Flanders 1960]. The ring $\mathcal{R}(w)$ in that case is also known as the Kronecker function ring.

*Corollary 7.4* *If a polynomial $p(X)$ belongs to an ideal $J$ of $\mathcal{R}(w)$, then all of the coefficients of $p(X)$ also belong to $J$.*

Proof. Because $w(p)$ is defined as the infimum of all the values $v(a_i)$ for coefficients $a_i$ of $p$, we have $w(a_i/p) \geq 0$ for each such coefficient, and hence $a_i = (a_i/p)p$ is an element of the ideal $J$. ■

*Corollary 7.5* *Let $A$ be a Prüfer or Dedekind domain and $v$ the product or divisor valuation on the field of fractions $K$ defined in terms of finitely-generated fractional ideals, and let $w$ be the extension of $v$ to the field $K(X)$ as in Theorem 3. There is a bijection between the ideals $J$ of $\mathcal{R}(w)$ and the ideals of $\mathcal{R}(v)$ given by $J \rightarrow J \cap \mathcal{R}(v)$.*

Proof.

(a) For any rational function $f$ in $K(X)$, observe that there is a polynomial $p$ in $K[X]$ such that $w(p) = w(f)$. Indeed, let $f = r/s$ where $r$ and $s$ are polynomials in $K[X]$. The fractional ideal $J$ generated by the coefficients of $s$ is invertible and represents the valuation $w(s)$ in the value group $G$ of finitely generated fractional ideals. Let $I$ be the inverse of $J$, then $I$ represents the valuation $w(t)$ of some polynomial $t$ (just take $t$ to be any polynomial where the coefficients of different monomials are a set of generators of $I$). We have $w(f) = w(r) - w(s) = w(r) + w(t) = w(r.t)$, and the polynomial $r.t$ has the same valuation as the rational function $f$.

Given any element $f(X)$ of an ideal $J$ in $\mathcal{R}(w)$, let $p(X)$ be a polynomial in $K[X]$ with the same valuation as $f(X)$. Because they have the same valuation, $p$ and $f$ generate the same principal ideal in $\mathcal{R}(w)$ and hence both are in $J$. Moreover, by corollary 3.4 above, all the coefficients of $p(X)$ are in $J \cap \mathcal{R}(v)$. Accordingly, $J$ is generated by $J \cap \mathcal{R}(v)$.

(b) We now show that for any ideal $\mathfrak{b}$ in $\mathcal{R}(v)$, we have $\mathfrak{b} = \mathfrak{b}\mathcal{R}(w) \cap \mathcal{R}(v)$.

Observe that the set $\mathfrak{b}$ has the following convex property: for any finite set of elements $x$, ..., $y$ in $\mathfrak{b}$, any element t in $\mathcal{R}(v)$ such that $v(\mathrm{t}) \geq \inf(v(x), ..., v(y))$ is also in $\mathfrak{b}$. In fact, take the fractional ideal of $A$ generated by the elements x,..., y of $\mathfrak{b}$. That fractional ideal represents $\inf(v(x), ..., v(y))$ in the value group $G$. Moreover, any element $t$ of $K$ such that $v(t) \geq \inf(v(x), ..., v(y))$ also belongs to that fractional ideal (recall that the ordering on the value group of $K$ is defined by reverse inclusion of the corresponding fractional ideals). Because that fractional ideal is contained in $\mathfrak{b}$, our observation follows.

Now consider any element $t$ in the ideal $\mathfrak{b}\mathcal{R}(w)$. We can write $t = px + ... + qy$, where $x$, ..., $y$ are elements of $\mathfrak{b}$, and $p$, ..., $q$ are elements of $\mathcal{R}(w)$. We have $w(t) \geq \inf(w(px), ..., w(qy)) \geq \inf(w(x), ..., w(y)) = \inf(v(x), ..., v(y))$. If the element $t$ is in the intersection $\mathfrak{b}\mathcal{R}(w) \cap \mathcal{R}(v)$, then by the convex property of the set $\mathfrak{b}$ noted above, the element $t$ also belongs to $\mathfrak{b}$. Therefore we see that $\mathfrak{b}\mathcal{R}(w) \cap \mathcal{R}(v)$ is contained in $\mathfrak{b}$. Because $\mathfrak{b}$ is obviously contained in $\mathfrak{b}\mathcal{R}(w) \cap \mathcal{R}(v)$, we have $\mathfrak{b}\mathcal{R}(w) \cap \mathcal{R}(v) = \mathfrak{b}$.

(c) Write $\boldsymbol{\alpha}$ for the mapping $\mathrm{J} \to \mathrm{J} \cap \mathcal{R}(v)$ and $\boldsymbol{\beta}$ for the mapping $\mathfrak{b} \to \mathfrak{b}\mathcal{R}(w)$, then (a) shows that $\boldsymbol{\beta}.\boldsymbol{\alpha}$ is the identity mapping on ideals $J$ of $\mathcal{R}(w)$, and (b) shows that $\boldsymbol{\alpha}.\boldsymbol{\beta}$ is the identity mapping on ideals $\mathfrak{b}$ of $\mathcal{R}(v)$. Hence $\boldsymbol{\alpha}$ and $\boldsymbol{\beta}$ are inverse bijections between ideals of $\mathcal{R}(w)$ and ideals of $\mathcal{R}(v)$. ■

NOTE: The above bijection respects multiplication of ideals, and so can be regarded as an isomorphism between the ideal groups of $\mathcal{R}(w)$ and $\mathcal{R}(v)$.

*Corollary 7.6* *Let A be a unique factorization domain. Then A[X] is also a unique factorization domain.*

Proof. Let $v$ be the natural $\mathbb{Z}$-divisor valuation on $K$ (the field of fractions of $A$) defined by prime factorization and let $w$ be the extension of $v$ to $K(X)$ as in Theorem 5.

Consider the natural $\mathbb{Z}$-divisor valuation $t$ on $K(X)$ defined by factorization into irreducible polynomials. Putting $w$ and $t$ together, we have a $\mathbb{Z}$-divisor valuation $u$ on $K(X)$ defined by $u(f) = (w(f), t(f))$. The valuation ring $\mathcal{R}(u)$ is $\mathcal{R}(w) \cap \mathcal{R}(t) = \mathcal{R}(w) \cap K[X] = A[X]$.

Note that the $\mathbb{Z}$-divisor valuations $w$ and $t$ are both surjective. Moreover, because the valuations $w$ and $t$ are independent of each other, the combined $\mathbb{Z}$-divisor valuation $u$ can readily be seen to be surjective. Consequently, the valuation ring $\mathcal{R}(u)$ = A[X] is a unique factorization domain. ■

*Corollary 7.7 (Dedekind's Prague Theorem) Let $\mathbf{a}_0, ..., \mathbf{a}_m$ be a set of (m + 1) elements in some field K of characteristic 0, and $\mathbf{b}_0, ..., \mathbf{b}_n$ be another set of (n + 1) elements in the same field.*

*In the subring $E = \mathbb{Z}[\mathbf{a}_0, ..., \mathbf{a}_m, \mathbf{b}_0, ..., \mathbf{b}_n]$, define the elements $\mathbf{c}_k = \sum \mathbf{a}_i\mathbf{b}_j$ where the sum runs over all pairs of indices i, j such that $i + j = k$ $(0 \leq k \leq m+n)$. For any valuation w defined on the field of fractions of E, we have the following inequality for each of the (m+1)(n+1) elements $\mathbf{a}_i\mathbf{b}_j$*

$$w(\mathbf{a}_i\mathbf{b}_j) \geq \inf w(\mathbf{c}_k) \quad (0 \leq k \leq m+n).$$

*In particular, $w(\mathbf{a}_i\mathbf{b}_j) \geq 0$ for any valuation w where $w(\mathbf{c}_k) \geq 0$ for all k, i.e., each element $\mathbf{a}_i\mathbf{b}_j$ is integral over the subring $\mathbb{Z}[\mathbf{c}_0, ..., \mathbf{c}_{m+n}]$.*

Proof. Let $L$ be the field of fractions of the integral domain $E$, that is $L = \mathbf{Q}(\mathbf{a}_0, ..., \mathbf{a}_m, \mathbf{b}_0, ..., \mathbf{b}_n)$. For any valuation $w$ on $L$, we can extend $w$ to the function field $L(X)$.

Let inf $w(\mathbf{a})$ be the infimum or greatest lower bound of the family $w(\mathbf{a}_i)$ as $i$ runs through the index set from 0 to $m$, and similarly for inf $w(\mathbf{b})$ and inf $w(\mathbf{c})$. By applying the Kronecker-Gauss lemma, we know that $\inf w(\mathbf{a}) + \inf w(\mathbf{b}) = \inf w(\mathbf{c})$.

Accordingly, $w(\mathbf{a}_i\mathbf{b}_j) = w(\mathbf{a}_i) + w(\mathbf{b}_j) \geq \inf w(\mathbf{a}) + \inf w(\mathbf{b}) = \inf w(\mathbf{c})$. ■

As another application, we can show that the ring of all integers in a number field is a Dedekind domain. Specifically, we will show the following characteristic property for a Noetherian domain to be Dedekind.

*<u>Proposition 8.</u> Let A be the ring of all integers in a number field K. Then for any non-zero ideal I of A, we can find an non-zero ideal J such that the product ideal IJ is a principal ideal. In other words, any non-zero ideal of A is invertible.*

<u>Proof</u>. $A$ is a Noetherian ring, so any ideal $I$ of $A$ is finitely generated. Let $a_1, a_2, \dots, a_m$ be a set of generators of the non-zero ideal $I$. Consider the polynomial

$$F(X) = a_1X + a_2X^2 + \dots + a_mX^m$$

For each homomorphism of the number field K into $\overline{\mathbf{Q}}$, we can transform the coefficients of $F$ and obtain a conjugate copy of $F$. (There are as many of these homomorphisms as the degree of K over $\mathbb{Q}$.) The product of all these conjugate copies of $F$ (one of which is $F$ itself) is a polynomial $H(X)$ whose coefficients must be in $\mathbb{Q}$ because they are invariant under any automorphism of $\overline{\mathbf{Q}}$. In fact these coefficients of $H$(X) are integers because they are integral over $\mathbb{Z}$ (being sum of products of algebraic integers), and the only numbers in $\mathbb{Q}$ integral over $\mathbb{Z}$ are the integers. Let $P$ be the greatest common divisor of these integral coefficients of $H(X)$.

Write $H(X) = F(X).G(X)$, where the polynomial $G(X)$ also has coefficients in the ring $A$ because these coefficients are integral over $\mathbb{Z}$ and moreover invariant under any automorphism of $\overline{\mathbf{Q}}$ over $K$. Let $G(X) = b_1X + b_2X^2 + \dots + b_nX^n$, and let $J$ be the integral ideal of $A$ generated by these coefficients.

We claim that the product ideal $IJ$, which is generated by the elements $a_ib_j$, is the same as the principal ideal $(P)$.

Corollary 7.7 tells us that for any valuation $w$ defined on $K$, we have $w(a_ib_j) \geq$ the greatest lower bound of the set $w(c_k)$, where the integers $c_k$ are the coefficients of the polynomial $H(X)$. In addition, $w(c_k) \geq w(P)$ for any $c_k$ because $P$ divides all the coefficients $c_k$ and $w$ must necessarily be $\geq 0$ on $\mathbb{Z}$ due to the defining properties of valuation. Therefore $w(a_i b_j) \geq w(P)$ for all $a_i b_j$ and for all valuations $w$ defined on $K$. According to Proposition 4, that condition is equivalent to $P$ being a divisor of $a_i b_j$ in the integral closure of $\mathbb{Z}$ in $K$, which is $A$. In other words, the principal ideal $(P)$ contains the product ideal $IJ$.

Because $P$ is the greatest common divisor of the integers $c_k$, we can express $P$ as a $\mathbb{Z}$-linear combination of the integers $c_k$. But note that each $c_k$ is in the ideal $IJ$ because it is a sum of the products $a_i b_j$. So the number $P$ itself is also in the product ideal $IJ$, and therefore we must have $IJ = (P)$. ■

## <u>REFERENCES</u>: